\documentclass[12pt]{article}
\usepackage{amsfonts}
\usepackage{amsmath}
\usepackage{amssymb}
\usepackage{mathrsfs}
\usepackage{theorem}
\usepackage{graphicx}
\usepackage[margin=1.2cm]{caption}

\newtheorem{theorem}{Theorem}

\theorembodyfont{\rmfamily}

\newcommand{\eat}[1]{}

\newcommand{\tz}{{\tilde{z}}}
\newcommand{\tw}{{\tilde{w}}}
\newcommand{\HH}{{\mathscr{H}}}

\newcommand{\rr}{{\mathbb{R}}}
\newcommand{\zz}{{\mathbb{Z}}}

\newcommand{\F}{\textup{\sf{F}}}

\title{Cusp transitivity in hyperbolic 3-manifolds}
\author{Roger Vogeler}

\date{}

\newcommand{\foorp}{{\unskip\nobreak\hfil\penalty50
 \hskip1em\vadjust{}\nobreak\hfil $\Box$
 \parfillskip=0pt \finalhyphendemerits=0 \par}}

\begin{document} 

\maketitle

\begin{abstract} 

Let $M$ be a cusped finite-volume hyperbolic three-manifold with isometry group $G$. Then $G$ induces a $k$-transitive action by permutation on the cusps of $M$ for some integer $k\ge 0$. Generically $G$ is trivial and $k=0$, but $k>0$ does occur in special cases. We show examples with $k=1,2,4$. An interesting question concerns the possible number of cusps for a fixed $k$. Our main result provides an answer for $k=2$ by constructing a family of manifolds having no upper bound on the number of cusps.
\end{abstract}






\section{Introduction} 
\label{Intro}
An action of a group $G$ on a set $S$ is called $k$-transitive if, for every choice of distinct elements $x_1, ... , x_k \in S$ and every choice of distinct targets $y_1, ... , y_k \in S$, there is an element $g\in G$ such that $g(x_i)=y_i$.  The term \emph{transitive} means 1-transitive; actions with $k>1$ are \emph{multiply transitive}.  The number of elements in $S$ is the \emph{degree}. Transitive actions are common (for example, every group acts transitively on itself by left multiplication), while multiply-transitive actions are relatively rare. The theory is well developed; see \cite{Dixon:1996uq}.

It is obvious that the isometry group of a complete finite-volume hyperbolic three-manifold induces a permutation action on the set of cusps. In this paper we call such a manifold \emph{k-transitive} if the induced action is $k$-transitive. Note that this definition is of the `inclusive hierarchy' type. For instance, a 3-transitive manifold is automatically 2-transitive, and possibly 4-transitive as well.

Kojima \cite{kojima1988isometry} shows by construction that every finite group $G$ occurs as the isometry group of some closed hyperbolic 3-manifold. At one stage, the construction involves a cusped manifold, with cusps labeled by the elements of $G$. The induced permutation action on the cusp set amounts to left-multiplication of the labels by elements of $G$. Since this  action is  transitive, the manifold at this stage is 1-transitive by our definition, and we see that 1-transitive manifolds are plentiful and can have any number of cusps. But  the action is also free, which implies that it is never 2-transitive. Hence one may begin to wonder, for $k\ge 2$, how common and how challenging to construct are the $k$-transitive manifolds?

In Section \ref{Examples} we give some examples of 1-, 2- and 4-transitive manifolds, all constructed as complements of hyperbolic links in $S^3$. In Section \ref{2-trans} we show a construction that leads to our main result:
 \begin{theorem} 
 \label{MainThm}
For the class of 2-transitive manifolds, there is no upper bound on the number of cusps.
 \end{theorem}
In Section \ref{Conclusion} we summarize known results and mention some open questions.

\section{Examples with $k=1,2,4$}
\label{Examples}
All the manifolds in this section are link complements in $S^3$. The constructions are easy, so we give only brief descriptions. The key idea is that, since the geometric structure of a hyperbolic link complement is unique, symmetries of the link are always realized by isometries of the complement. Where no outside reference is given, hyperbolicity has been verified using SnapPea \cite{Weeks:1999uq}.

\begin{figure}[h]
\begin{center}
\includegraphics[scale=0.6]{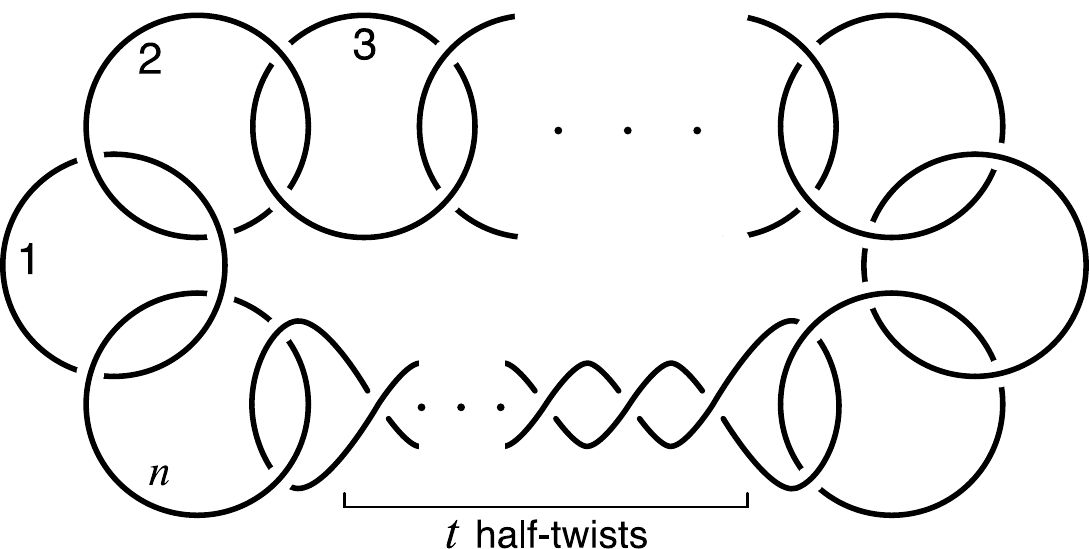}
\end{center}
\caption{Chain link with $n$ components and $t$ half-twists.}
\label{ChainLink}
\end{figure}

\subsection{A chain construction}
Consider a cycle of $n$ unknotted loops, $n\in\zz ^+$, linked sequentially, and having $t$ half-twists, $t\in\zz$, as in Figure \ref{ChainLink}.  The chain itself, viewed as a single loop, must be unknotted. There is an obvious 1-transitive cyclic action on the link components. 

Two sub-families of such chains appear as examples in Thurston's notes, and the complete family is analyzed by Neumann and Reid \cite{neumann1992arithmetic}. They show that for $n\ge5$ the link is hyperbolic for every value of $t$, and for $n < 5$ the link is hyperbolic for all but $5-n$ values of $t$, the exceptions being the cases with the least  overall twisting.  

\begin{figure}[h]
\begin{center}
\includegraphics[scale=0.5]{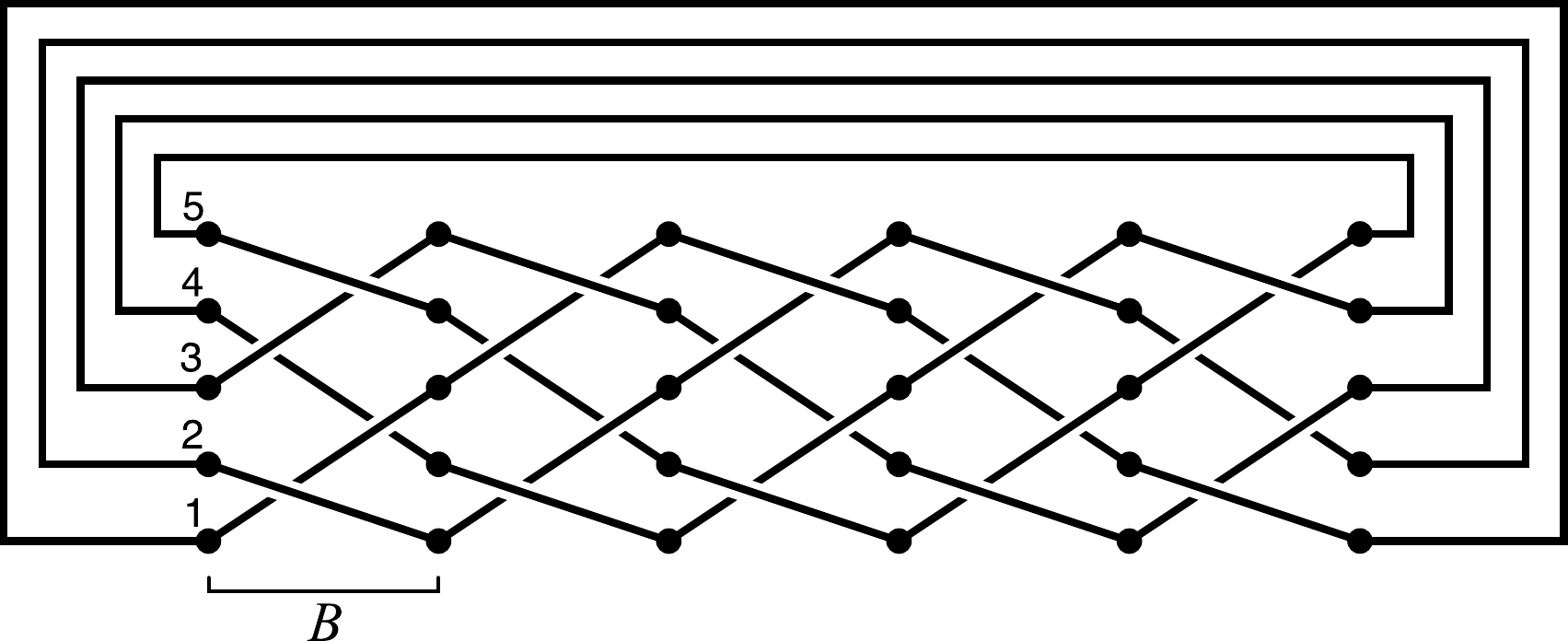}
\end{center}
\caption{Braid $B$ has associated permutation $(1 \,3 \,5 \,4 \,2)$. Although the closure of $B$ is the unknot, the closure of $B^5$ is hyperbolic. Shifting the diagram cyclically one step to the right gives a 1-transitive action on the five components.}
\label{B5closure}
\end{figure}

\subsection{A braid construction}
Suppose $B$ is a pseudo-Anosov braid on $n$ strands for which the associated permutation is cyclic of order $n$. Then the braid $B^n$ induces the trivial permutation on the strands, and the closure of $B^n$ is thus a link having $n$ components. Again, there is an obvious 1-transitive cyclic action. An example is shown in Figure \ref{B5closure}.

The closed braid might fail to be hyperbolic. In particular, since a braid is closed by performing a $(1,0)$ Dehn filling on the braid axis, it might happen that the meridian is too short to allow a hyperbolic filling. In this case, the longer braid $(B^n)^m$ will have meridian $m$ times as long, and thus, by the `$2\pi$ theorem' of Gromov and Thurston \cite{bleiler1996spherical}, the closed braid will be hyperbolic for sufficiently large $m$. (It seems likely that the $n$th power of $B$ is already enough to guarantee a meridian longer than $2\pi$, but there may be exceptions to this.)

\subsection{Two cubical constructions}

\begin{figure}[h]
\begin{center}
\includegraphics[scale=0.25]{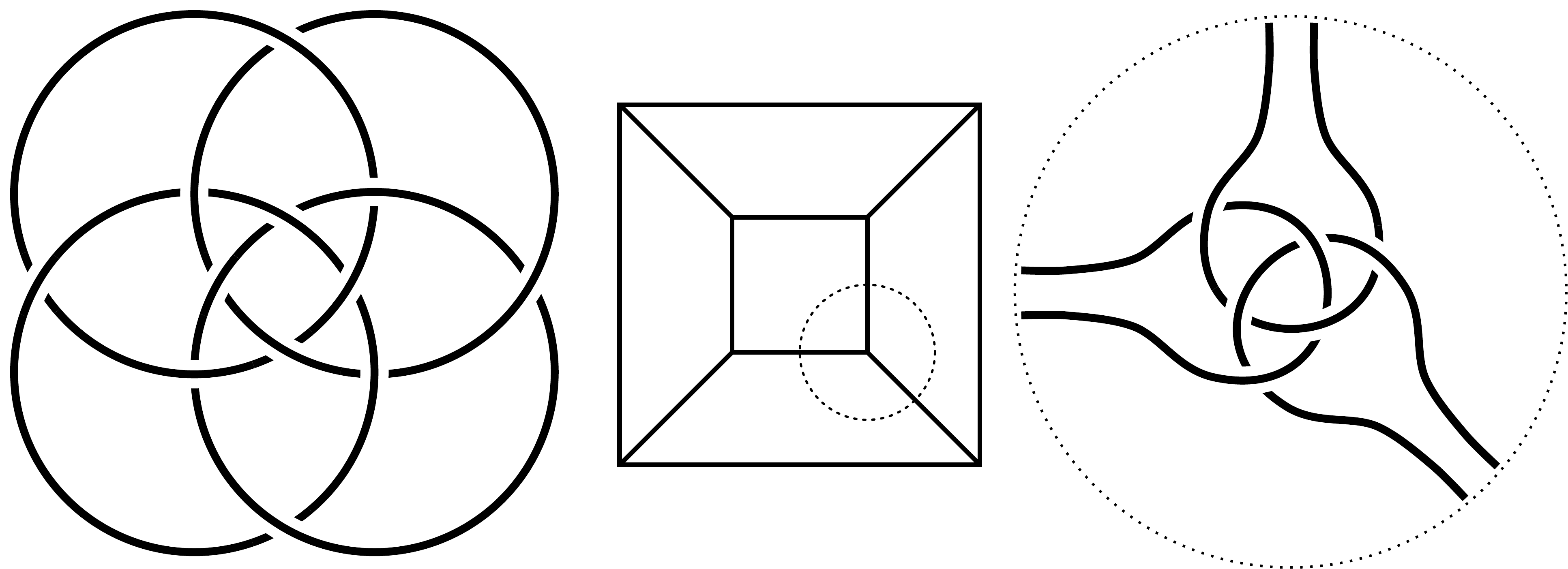}
\end{center}
\caption{Left: The complement of this link is  4-transitive. Center and right: Each vertex of a cube is replaced by three interlocking half-loops. The complement is 1-transitive, with 12 cusps.}
\label{CubeLink}
\end{figure}

Consider the four planes in $\rr ^3$ given by $Ax + By + z=0$, with $A,B\in\{\pm 1\}$. These planes are perpendicular to the four diagonals of the cube with vertices $(\pm 1, \pm 1, \pm 1)$. Intersecting these planes with the unit sphere gives four great circles  intersecting at twelve points. By resolving each intersection to an over- or under-crossing, in alternating fashion, a hyperbolic link is obtained. A planar diagram is shown in Figure \ref{CubeLink} (Left).

The group of orientation-preserving symmetries of the cube carries over naturally to the link components, and it is an easy exercise to check that this action is 4-transitive. This link can also be obtained by the braid construction given above, but the 4-transitivity is not so obvious.

An analogous construction, where the cube is replaced by an icosahedron, yields a 2-transitive manifold having 6 cusps.

The second cubical construction builds a 12-component link by following the structure of the edges of a cube, as shown in Figure \ref{CubeLink}. The complement is hyperbolic and 1-transitive. In contrast to the chain and braid examples shown earlier, the action is non-cyclic. 

\eat{ 
\subsection{A conjecture}
How many cusps can a $k$-transitive manifold have? The examples shown in Section \ref{} strongly suggest that any positive number of cusps can occur in a 1-transitive manifold. For $k=2$ the situation seems more difficult, but the construction given in Section \ref{2-trans} will show that the possibilities include every prime power. The situation feels quite different for $k\ge 3$, due in part to the scarcity of groups that can act in the required way, and in part to the spatial restrictions inherent in dimension three. While hoping that it is not the case, we suspect that there is in fact an upper bound on the number of cusps whenever $k\ge 3$. 
} 

\section{A family of 2-transitive manifolds}
\label{2-trans}
In this section we construct an infinite family of 2-transitive manifolds for which the number of cusps increases without bound. The construction builds on a family of regular maps described by Biggs \cite{Biggs:1971kx} and further studied by James and Jones \cite{James:1985uq}.


Here is an overview of the construction. Let $S$ be a surface carrying one of the regular maps of Biggs. There is a natural 2-transitive action on the $n$ faces of this map by its automorphism group. The product $S\times S^1$ is formed, and within this non-hyperbolic manifold a system of $n$ simple closed curves is specified. The curves are carefully designed to maintain the original 2-transitive symmetry. The complement of this link is the desired manifold $M$. Finally, to verify that $M$ is hyperbolic, it is viewed as a punctured-surface bundle over $S^1$ with pseudo-Anosov monodromy.

\begin{figure}[h]
\begin{center}
\includegraphics[scale=0.63]{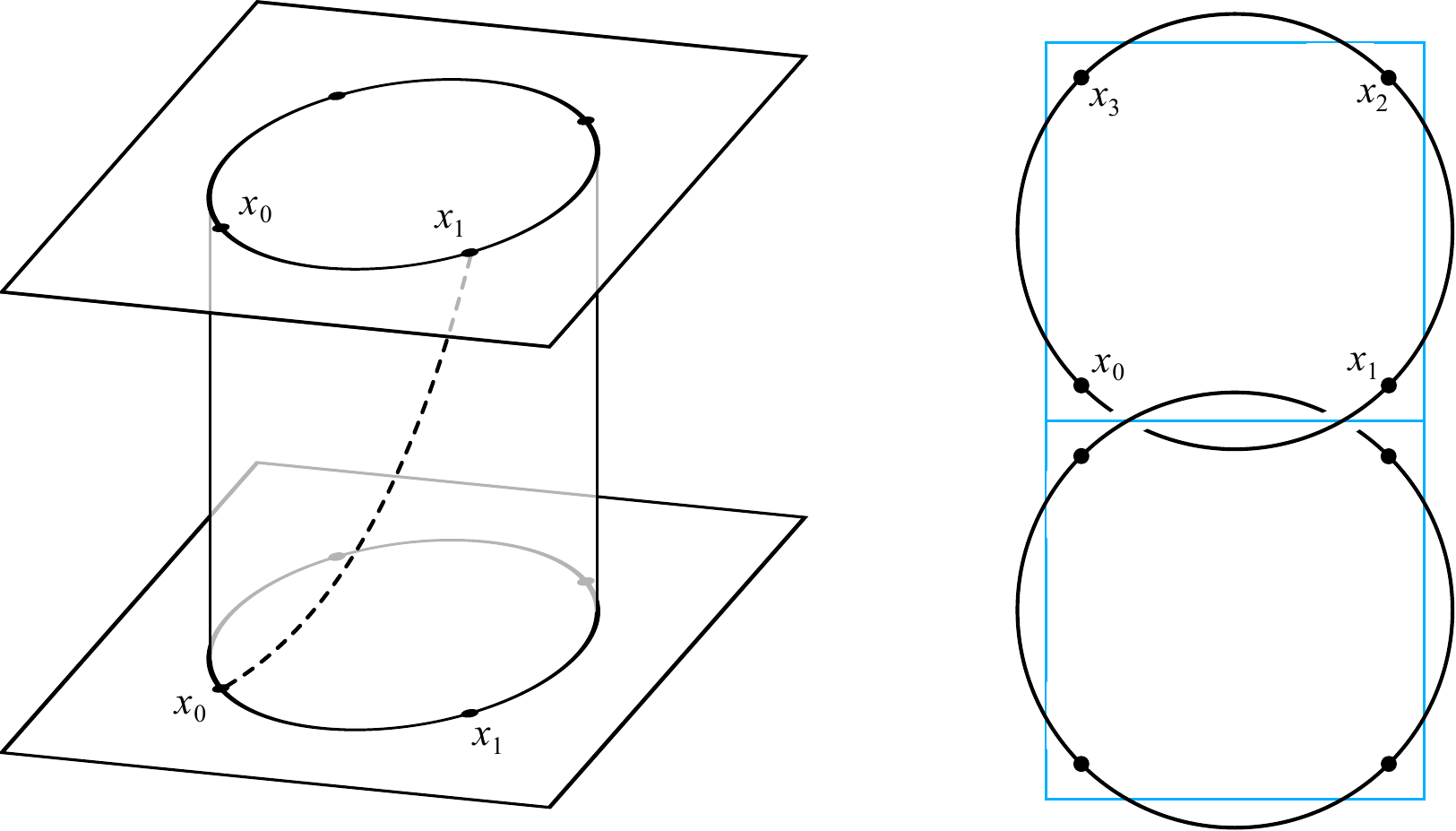}
\end{center}
\caption{Left: Each arc is a helical segment joining the image of $x_i$ in the lower boundary to the image of $x_{i+1}$ in the upper boundary. After boundary identification, four arcs become a single closed curve. Right: The view from above, after increasing the radius so that each arc extends out to link with an arc of the neighboring face. }
\label{HelicalArc}
\end{figure}

We now describe the construction in more detail. Let $n$ be a prime power greater than 3, and let $\F$ be the finite field of order $n$. Biggs uses the arithmetic of $\F$ to construct orientable regular maps (tessellations) on a closed surface.  Such a map has $n$ faces, each face a $(n-1)$-gon. The genus of the surface is given by 
\begin{equation*}
g_n=\left\{
     \begin{array}{ll}
       1+\frac{\textstyle n(n-7)}{\textstyle 4} & \textrm{if} \;\; n\equiv 3 \mod 4\\
       1+ \frac{\textstyle n(n-5)}{\textstyle 4} & \textrm{otherwise.}
     \end{array}
    \right .
\label{genus}
\end{equation*}
The faces admit a natural labeling by the elements of $\F$, and the automorphism group of the map is exactly the natural action of the affine group $A(1,\F)$ on the face labels. The 2-transitivity of this action is an easy exercise. The significance of this, combinatorially, is that each face of the map shares an edge with each of  the remaining faces.

We take $S$ to be the orientable surface carrying this regular map, endowed with the natural constant-curvature metric (see \cite{edmonds1982regular}). The case $n=5$ is used in the illustrations; the surface is a square torus tiled by five square faces. 

The construction continues inside a preliminary three-manifold $M_I=S\times I$, where $I$ is the unit interval. For each face $f$ of $S$, $n-1$ arcs are specified in $M_I$ as follows. First, let $P$ be the center point of $f$, and let $C$ be the circle of radius $\rho$ and circumference $\sigma$ in $S$ centered at $P$. Let the vertices of $f$ be denoted $v_i$ sequentially around the perimeter, for $0\le i\le n-2$. Let $x_i$ be the point of intersection of $C$ and the radial segment extending from $P$ to $v_i$.  Now $C_I=C\times I$ is a cylinder surrounding the fiber $\{P\}\times I$ in $M_I$. The desired arcs are defined to be the helical segments $s_i$ of slope $\frac{n-1}{\sigma}$ lying on $C_I$, each extending from the image of $x_i$ in $f\times\{0\}$ to the image of $x_{i+1}$ in $f\times\{1\}$ (indices modulo $n-1$). Figure \ref{HelicalArc} illustrates this step.

\begin{figure}[h]
\begin{center}
\includegraphics[scale=0.8]{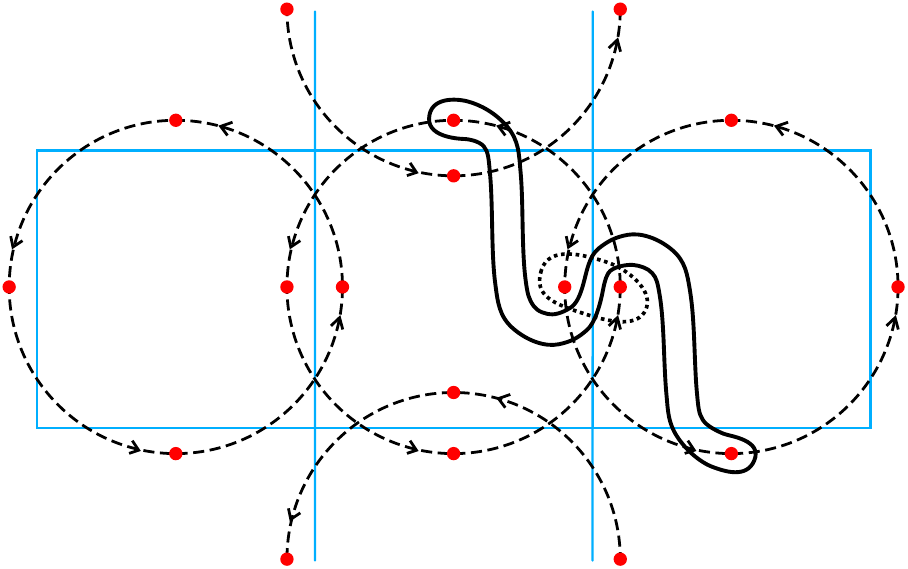}
\end{center}
\caption{The point-pushing homeomorphism $\mathscr{H} $. The marked points are punctures. Dashed arrows represent the point-pushing trajectories. The small dotted loop encloses two punctures; the long snaking loop is its image under $\mathscr{H} $. With further iterations it would limit on the $f_+$ foliation in the sense of Thurston \cite{thurston1988geometry}.}
\label{PointPush}
\end{figure}

Next, $M_\circ$ is formed by identifying the two boundary components of $M_I$ by equality in the first factor; hence $M_\circ$ is the product $S\times S^1$. By this identification the cylinder $C_I$ becomes a torus, and the segments $s_i$ join end-to-end yielding a closed curve, namely, a $(n-1,1)$ torus knot.

Finally, to obtain the desired link, we adjust the size of circle $C$ as follows. A face of $S$ has two characteristic radii:  $r_1$, the distance  from center $P$ to an edge midpoint, and $r_2$, the distance from $P$ to a vertex. By requiring the radius $\rho$ to satisfy $r_1 < \rho < r_2$, each copy of $C$ lies partly on its own face and partly on the neighboring faces. Hence each strand of the corresponding torus knot extends radially outward just far enough to link with a strand of the neighboring knot. Manifold $M$ is now defined to be the complement in  $M_\circ$ of this link.


\eat{
It is easy to see that $M$ is the mapping torus of a point-pushing homeomorphism $\mathscr{H}$ of a punctured surface. The fiber is just the original surface $S$, but with $n(n-1)$ points missing where the surface intersects the strands of the link. We view $\mathscr{H} $ as acting on the fiber  at level $\frac{1}{2}$ with respect to the original $I$ factor, as shown in Figure \ref{PointPush}.
}

In $M$, each original cross section $S\times \{x\}$ (for $x\in S^1$) is now a punctured surface, since the $n(n-1)$ points of intersection with the strands of the link are now missing. For convenience in the subsequent analysis, we specify the cross section at level $\frac 1 2$ with respect to the original $I$ factor as the `base surface' and observe that $M$ is the mapping torus of a point-pushing homeomorphism $\HH$ of this surface. A diagram for $\HH$ is shown in Figure \ref{PointPush}.

\begin{figure}[h]
\begin{center}
\includegraphics[scale=0.38]{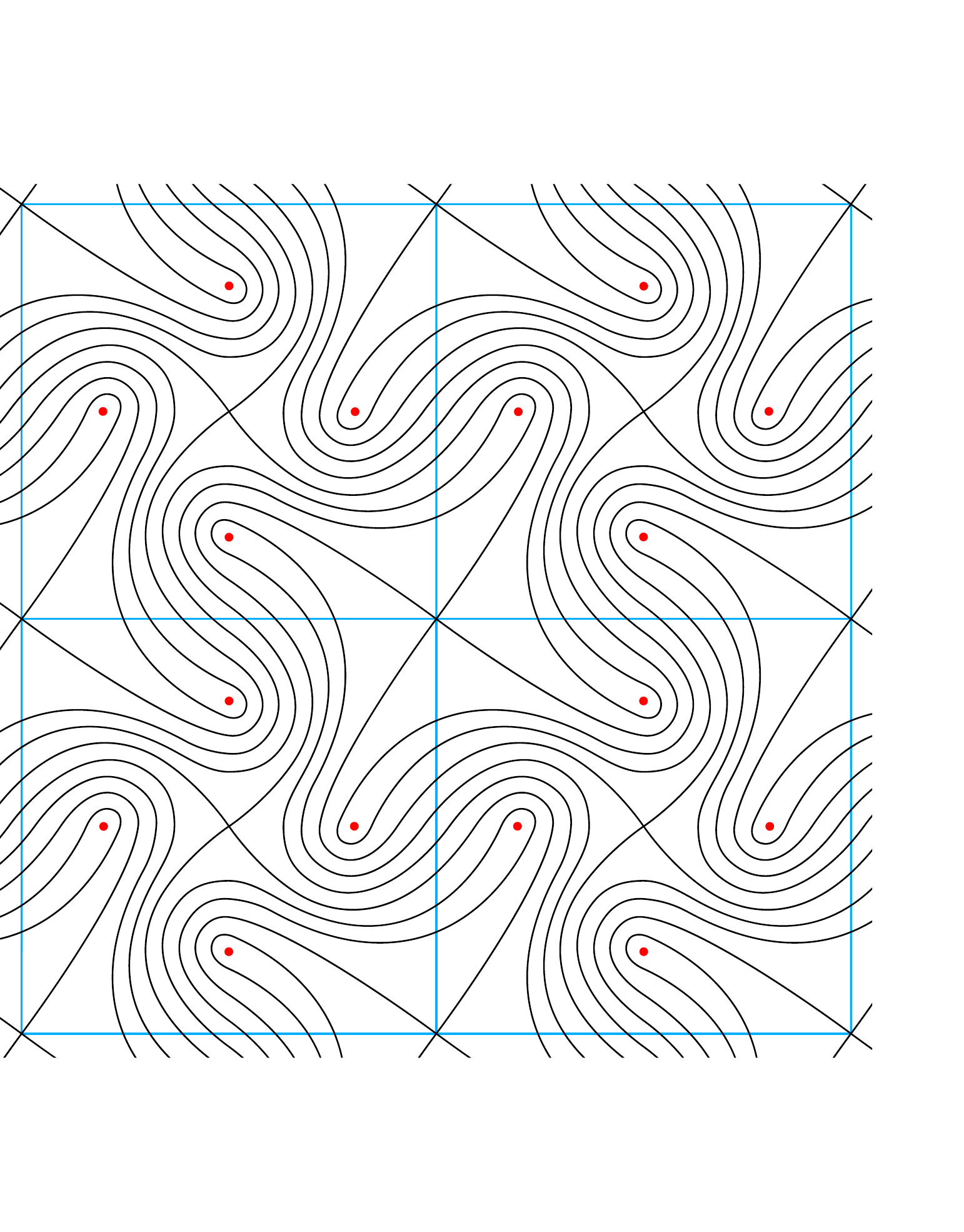}
\includegraphics[scale=0.38]{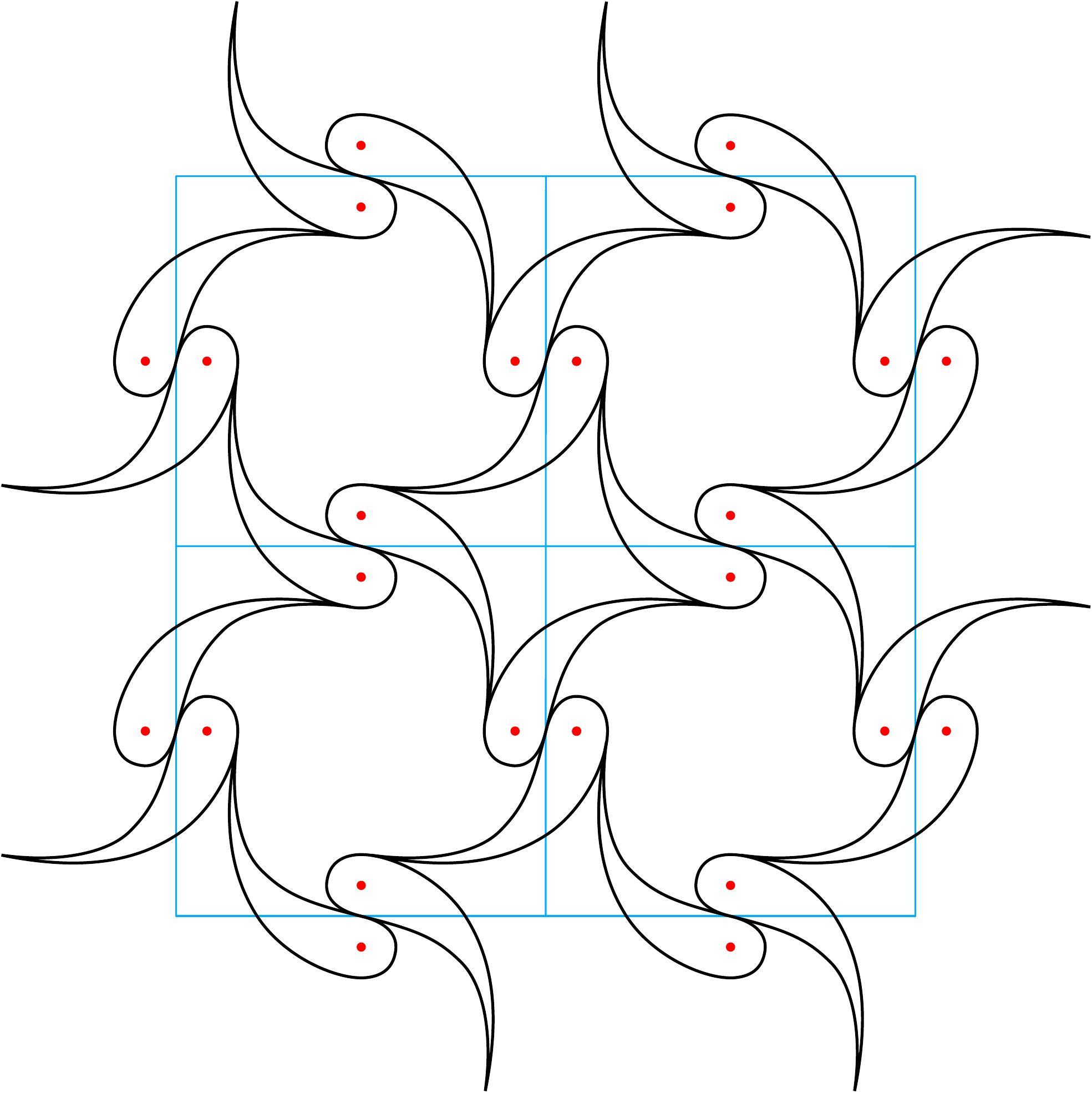}
\end{center}
\caption{Left: An approximation of the $f_+$ foliation, showing the 4-pronged singularities. Right: A train track for $\mathscr{H}$.}
\label{TrainTrack}
\end{figure}

A homeomorphism is pseudo-Anosov if there exists a pair of  transverse singular foliations, $f_+$ and $f_-$, invariant under $\mathscr{H} $, and an associated constant $\lambda>1$, the dilatation. The surface is stretched and compressed along the two foliations by factors of $\lambda$ and $\lambda ^{-1}$, respectively. By a result of Thurston \cite{Thurston:1982fk}, the mapping torus of a pseudo-Anosov homeomorphism is hyperbolic.

\eat{
Figure \ref{TrainTrack} shows an approximation of the $f_+$ foliation for $\HH$. Singularities occur at the vertices and face-centers of the underlying map; these points are fixed by $\mathscr{H}$. A train track for $\mathscr{H}$, as discussed by Bestvina and Handel \cite{bestvina1995train}, 
is also shown. Positive weights for the track will be computed below, along with $\lambda$. The complement of the track consists of once-punctured monogons and cusp-cornered polygons. It follows that $\mathscr{H}$ is indeed pseudo-Anosov. Details for the theory of train tracks are given by Penner and Harer \cite{Penner:1992vn}.
}

Figure \ref{TrainTrack} shows an approximation of the $f_+$ foliation for $\HH$. Singularities occur at the vertices and face-centers of the underlying map; these points are fixed by $\HH$. A train track for $\HH$, as discussed by Bestvina and Handel \cite{bestvina1995train}, 
is also shown.  The structure of the track is equivariant with respect to the symmetries of the surface. The complement of the track consists of once-punctured monogons and cusp-cornered polygons. The track captures in a discrete way the overall structure of the $f_+$ foliation, and is mapped, with appropriate stretching, onto itself by $\HH$ (modulo some adjustment homotopic to the identity). Details for the theory of train tracks are given by Penner and Harer \cite{Penner:1992vn}.

We turn next to the calculation of $\lambda$ and the weights associated to branches of the track. A few different approaches to this problem are possible. The underlying symmetry allows some simplifications to be made. 

\begin{figure}[h]
\begin{center}
\includegraphics[scale=0.7]{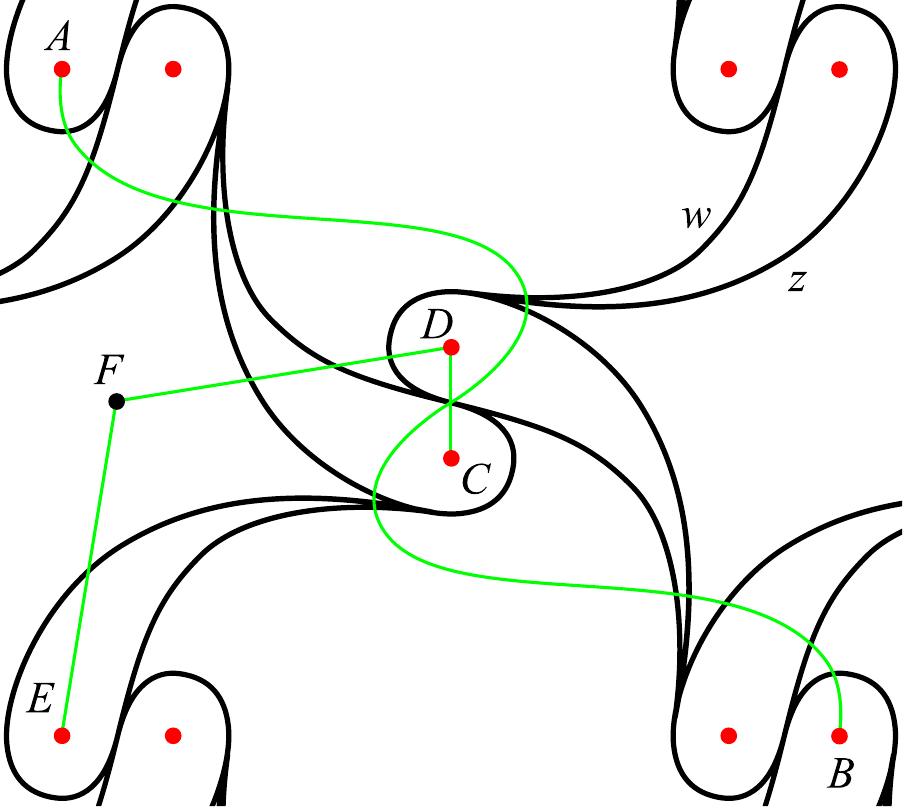}
\end{center}
\caption{Train track and arcs for computing the dilatation $\lambda$.}
\label{Dilatation}
\end{figure}

\eat{
To compute the track weights and the dilatation $\lambda$ for $\mathscr{H}$ we refer to the diagram in Figure \ref{Dilatation}. Points labeled $A$ through $E$ are punctures, while $F$ is a fixed point of $\mathscr{H}$. Each branch of the track carries a weight. Two branches are assigned weights $w$ and $z$; the weights for all other branches are determined by symmetry and by the summation rule where branches meet. The weights correspond to a measure for transverse arcs. 
}

In Figure \ref{Dilatation} the points labeled $A$ through $E$ are punctures, while $F$ is a fixed point of $\HH$. Each branch of the track is assigned a positive weight. The weights determine a crossing measure for transverse arcs. Two branches in the diagram are assigned weights $w$ and $z$; the weights for all other branches are determined by symmetry and by the summation rule where branches merge. Thus each semicircular branch that half-surrounds a puncture has weight $w+2z$ and each very short branch found between a pair of nearest-neighbor punctures has weight $2w+2z$.

Arc $\overline{AB}$ crosses 17 branches and has total crossing measure $10w+14z$. Its image under $\mathscr{H}$ is arc $\overline{CD}$, which has crossing measure $2w+2z$. One may view this as compression by $\HH$ in the transverse ($f_-$) direction, resulting in the full weight $10w+14z$ being forced onto a portion of  track that previously carried just $2w+2z$. Since the contraction factor for the  compression is $\lambda^{-1}$, one may attribute the increased `density' to a scaling by $\lambda$; 
hence
\begin{equation}\label{weights1}10w+14z=\lambda (2w+2z).\end{equation}
Similarly, arc $\overline{DF}$ has measure $2w+3z$, while its image, arc $\overline{EF}$, has measure $z$, and we obtain 
\begin{equation}\label{weights2}2w+3z=\lambda z.\end{equation}
With these two equations, we have enough constraints to determine the unknowns. Since weights are determined only up to a scale factor, we can set $z=1$ and solve, first obtaining $w=\sqrt{2}$, and then 
$$\lambda =3+2\sqrt{2}$$
and
$$ \lambda^{-1}=3-2\sqrt{2}.$$

This quick solution may seem a bit ad hoc; for a more systematic view of things let us divide both sides of equation \eqref{weights1} by 2 and then subtract equation \eqref{weights2} to produce
\begin{equation}\label{weights3} 3w+4z=\lambda w. \end{equation}
Equations \eqref{weights2} and \eqref{weights3} together now form a linear system 
$$\left[\begin{array}{rr}
3 & 4 \\
2 & 3 \\
\end{array}\right]
\left[\begin{array}{r}
w \\
z \\
\end{array}\right]
= \lambda
\left[\begin{array}{r}
w \\
z \\
\end{array}\right],$$
and we see how the weights constitute an eigenvector for  eigenvalue $\lambda$.

A related (and the most standard) approach is based on viewing the train track as an embedded graph and analyzing the self-map induced on it by $\HH$. Every branch is assigned a numerical value; these correspond to a tangential measure, and may be interpreted simply as edge lengths. The map stretches each edge onto some sequence of edges; this data is encoded by an incidence matrix. The stretched length is then some particular sum of original edge lengths, and one wishes to attribute all the stretching to a single scale factor. This, of course, is again an eigenvalue problem.

\begin{figure}
\begin{center}
\includegraphics[scale=0.25]{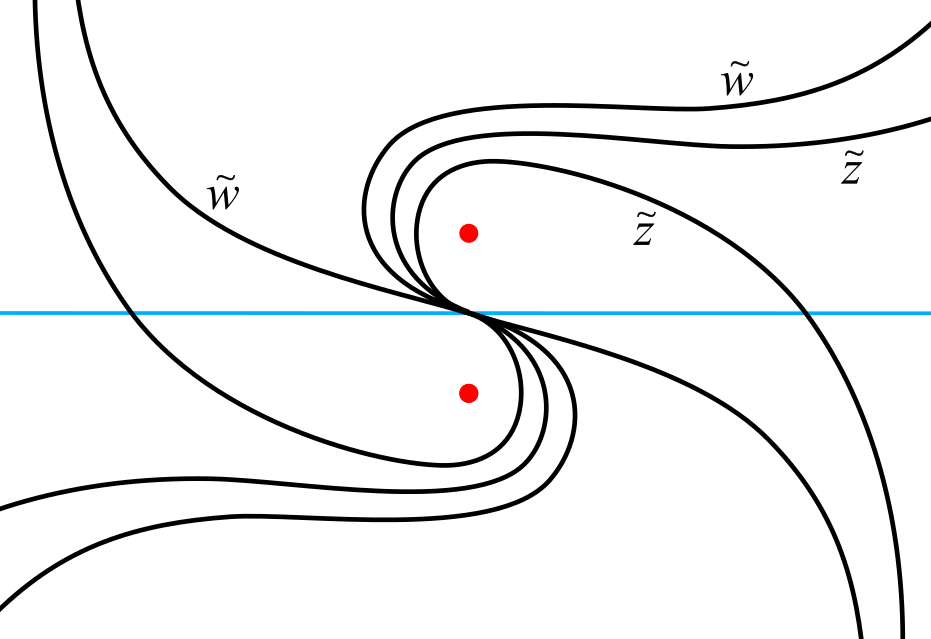}
\end{center}
\caption{The modified track has just two classes of branches. The semicircular branches have been completely split, and each  short branch has been reduced to a switch located at the edge-midpoint of a face. }
\label{TrackMod}
\end{figure}

The full incidence matrix of a train track for $\HH$ is rather large, but two easy simplifications make the problem tractable. First we modify the track as shown in Figure \ref{TrackMod}, so there are just two classes of branches. These are labeled $\tw$ and $\tz$. Then we replace the incidence matrix by a two-by-two matrix that imposes the same constraints on edge lengths that are implicit in the incidence matrix. Using the diagrams in Figures \ref{PointPush} through \ref{TrackMod} it is rather straightforward to trace out the image of a $\tw$ edge: it is stretched onto a sequence of edges of type $\tw,\tz,\tw,\tz,\tw$. Similarly, a $\tz$ edge is mapped onto $\tw,\tz,\tw,\tz,\tw,\tz,\tw$. Hence we have 
$$3\tw + 2\tz = \lambda\tw \quad\quad \textup{and}\quad\quad 4\tw + 3\tz = \lambda\tz.$$
The matrix for this linear system is just the transpose of the one that occurred earlier, so the eigenvalues are unchanged. The weights $\tw$ and $\tz$ for the tangential measure are the reciprocals (up to a scale factor) of those obtained earlier for the transverse measure. An overview of these concepts, placing them in the larger context of Thurston's work, can be found in the Epilogue of \cite{Penner:1992vn}.

Yet another approach to analyzing $\HH$ has been generously pointed out by an anonymous referee. Due to the large symmetry group $A(1,\F)$, the associated quotient space of $S$ is quite simple. A fundamental domain consists of two neighboring triangles in the barycentric subdivision of a face, and the induced edge identifications yield a topological sphere with three cone points, something like a triangular pillow or a turnover pastry. The punctures all get identified to a single point of the quotient; removing this point, along with the three cone points, produces a four-punctured sphere. This surface has a natural metric coming from its two-covering by a flat torus via the classical elliptic involution, with the four fixed points corresponding to the punctures. The image of $\HH$ in this reduced setting can be analyzed explicitly using Dehn twists, and its lift to the torus is a genuine Anosov transformation having scale factors $3\pm2\sqrt{2}$. Details for these steps can be found in \cite{bestvina1995train}, \cite{thurston1988geometry}, and the  primer \cite{farb2011primer} by Farb and Margalit.

The details used in the above approaches (specifically, the diagrams in Figures \ref{Dilatation} and \ref{TrackMod}, and the nature of the quotient space) do not depend on the faces being squares. Thus the computations work out as shown not just for $n=5$, but for every prime-power value of $n$. Hence the same dilatation occurs in every case.

\section{Further questions}
\label{Conclusion}
The constructions shown in Section \ref{Examples}, along with Kojima's method, suggest that 1-transitive manifolds may allow for quite a diversity of structures. It would be interesting to find a more systematic understanding of the big picture.

For $k=2$, more examples can be produced by easy variations of the above construction based on the Biggs surfaces. Are there other 2-transitive manifolds having substantially different structure? Is every 2-transitive group action realized by some 2-transitive manifold?

There are lots of 3-transitive manifolds having just three cusps, but this is not too surprising since a 3-transitive action on a set of three objects is not really very exotic. In particular, the three-component chain links described earlier all have obvious 3-transitive symmetry, and all but two of them are hyperbolic. The famous Borromean rings are another example. For more than three cusps, the situation appears much harder. The author is aware of one 3-transitive manifold with six cusps, and another with eight, but no others. 

 For $k=4$, our only known example is the manifold described in Section \ref{Examples}, and easy variations of it, all having four cusps. 
 
 \eat{
 For $k=5$, not a single example is known to the author. 
 
 It seems reasonable to conjecture that there is a largest $k$ for which $k$-transitive manifolds exist, and that for each $k\ge 3$ there is an upper bound on the possible number of cusps.

In the original version of this paper we posed the following conjecture which has since been answered in .... We include the original statement of the conjecture to accurately reflect the history of the problem. "
}

In the original version of this paper we stated two conjectures: that there is a largest $k$ for which $k$-transitive manifolds exist, and that for each $k\ge 3$ there is an upper bound on the possible number of cusps. Both have since been proven, and sharp bounds provided, by Ratcliffe and Tschantz \cite{ratcliffe2019cusp}.

\bibliographystyle{amsplain}		
\bibliography{CuspTransRefs}		

\smallskip


\smallskip
\noindent{\small\textsc{Department of Mathematical Sciences, 
Central Connecticut State University, 
New Britain, Connecticut 06050, USA}}

\smallskip
\noindent{\small\tt{vogelerrov@ccsu.edu}}

\end{document}